\documentclass[12pts,reqno]{amsart}
\usepackage[margin=1in]{geometry}
\usepackage{amssymb, amscd, amsmath, url}
\usepackage{lipsum}
\usepackage{etoolbox}
\usepackage[numbers]{natbib}
\usepackage{natbib}
\patchcmd{\subsection}{-.5em}{.5em}{}{} 
\patchcmd{\section}{\normalfont}{\bfseries\Large}{}{}

\numberwithin{equation}{section}
\usepackage{listings}
\usepackage{matlab-prettifier}      
\usepackage{graphicx}
\usepackage[foot]{amsaddr}
\usepackage{float}
\usepackage{graphicx}
\usepackage{lmodern}

\begin{document}
	
	\title[]{Advancements in Spectral Collocation Methods for High-Order Eigenvalue Problems}
	
	
	
	\author{\bfseries Sameh Gana$^*$ }
	\address{Department of Basic Sciences, 
		Deanship of Preparatory Year and Supporting Studies, 
		Imam Abdulrahman Bin Faisal University, P.O. Box 1982, Dammam, 34212, Saudi Arabia.}
	\address{$^*$Corresponding author:\textnormal{sbgana@iau.edu.sa}}
	
	\begin{abstract}
		
 This paper focuses on computing spectral solutions for high-order eigenvalue problems using an efficient discretization method based on Chebfun spectral discretization algorithms and domain truncation. We solve several numerical eigenvalue problems, demonstrating both the accuracy and computational efficiency of the proposed approach.

	\end{abstract}
	\maketitle
	
	\section{Introduction}
	 
Eigenvalue problems arise in many areas of physics, ranging from solving classical electromagnetic problems to calculating the quantum bound states of the hydrogen atom. In quantum mechanics, the solutions of the radial Schrödinger equation describe the bound-state energies of the non-relativistic hydrogen atom. Eigenvalue problems commonly appear in mechanical vibration analysis and stability studies in hydrodynamics and magnetohydrodynamics. Due to their complexity, numerical methods are often necessary to obtain accurate spectral solutions. For more applications, see ~\cite{Orszag, Gana1,Gana3,Gana5, Mhadhbi1,Mhadhbi2}.

In recent years, there has been significant interest in developing the classical spectral collocation method, which has become an effective numerical tool for eigenvalue problems.  Spectral methods provide exponential convergence for several problems, generally with smooth solutions. However, these numerical methods are not suitable for approximating high-index eigenvalues, particularly in the context of singular eigenvalue problems.

Modern advances in computational mathematics have led to the development of specialized software tools, including MATSLISE ~\cite{Ledoux}, SLEDGE ~\cite{Fulton}, and SLEIGN ~\cite{Bailey}. Despite their improvements, these conventional methods remain inadequate in computing eigenfunctions. Their limitations are evident in the presence of singularities or high-order problems.

Integrating Chebfun algorithms with spectral methods provides greater flexibility for solving a range of eigenvalue problems. 
Additionally, the software has been successfully applied to problems involving singularities, where classical numerical methods often struggle  ~\cite{Driscoll}.
For instance, Driscoll and Trefethen ~\cite{TrefethenS} demonstrated the effectiveness of Chebfun  in achieving high-precision eigenvalue solutions for high-order differential problems. For more complete descriptions of the Chebyshev Collocation method and more details on the Chebfun software system, we refer to
~\cite {TrefethenS, CANUTO,FORNBERG,TrefethenF,TrefethenT,Gheorghiu1,Gheorghiu2,Gana2,Gana4}.

The purpose of this paper is to assert that Chebfun, along with the spectral collocation methods,
can provide accuracy, robustness and simplicity of implementation. In addition, these methods can compute the whole set of eigenvectors and provide some details on the accuracy and numerical stability of
the results provided. 
 
The structure of this paper is organized as follows. Section 2 reviews the fundamental concepts underlying the Chebfun system and the Chebyshev spectral collocation method. Section 3 is the core part of our study. We analyze a series of numerical experiments involving problems with a mixture of boundary conditions. These numerical computations will involve high-order eigenvalue problems, allowing us to evaluate the performance of the spectral collocation method and Chebfun across various scenarios. We conclude with Section 4, which is devoted to conclusions and open problems.

	\section{Chebfun System and Chebyshev Spectral Collocation Methodology}
	
	The Chebfun system contains algorithms that yield high-precision resolution spectral collocation methods on Chebyshev grids. The Chebops tools in the Chebfun system for solving differential equations are summarized in ~\cite {Driscoll} and ~\cite {Hale}. \\
The Chebops implementation effectively merges the numerical analysis approach of spectral collocation with the computational methodology of spectral discretization matrices. The Chebfun system explained in ~\cite {TrefethenS} solves the eigenproblem by choosing a reference eigenvalue and checks the convergence of the process.

	The Spectral Collocation method for solving differential equations consists of constructing weighted interpolants of the form:
	
	\begin{equation}\label{eq1}
		u(x)\approx P_{N}(x)=\sum_{j=0}^{N}\frac{\alpha(x)}{\alpha(x_{j})}\phi_{j}(x)u_{j},
	\end{equation}	
	where $x_{j}$ for $j=0,....,N$ are interpolation nodes, $\alpha(x)$ is a weight function, $$u_{j}=u(x_{j}),$$ and the interpolating functions
	$\phi_{j}(x)$ satisfy $$\phi_{j}(x_{k})=\delta_{j,k}$$ and $$u(x_{k})=P_{N}(x_{k})$$ for $k= 0,....,N$.\\ Hence, $P_{N}(x)$ is an interpolant of the function $u(x)$.\\ By taking $l$ derivatives of \ref{eq1} and evaluating the result at the nodes  $x_{j}$, we get:
	$$u^{(l)}(x_{k})\approx \sum_{j=0}^{N}\frac{d^{l}}{dx^{l}}\left[\frac{\alpha(x)}{\alpha(x_{j})}\phi_{j}(x)\right]_{x=x_{k}}, \quad k= 0,....,N.$$
	The entries define the differentiation matrix:  
	$$D^{(l)}_{k,j}= \frac{d^{l}}{dx^{l}}\left[\frac{\alpha(x)}{\alpha(x_{j})}\phi_{j}(x)\right]_{x=x_{k}}.$$
	The derivatives values $u^{(l)}$ are approximated  at the nodes $x_{k }$ by $ D^{(l)}u $.\\
	The derivatives are converted to a differentiation matrix form and the differential equation problem is transformed into a matrix eigenvalue problem.
	
	The eigenfunctions 	$u(x)$ of the eigenvalue problem approximate finite terms of Chebyshev polynomials as
	\begin{equation} \label{eq2}P_{N}(x)= \sum_{j=0}^{N}\phi_{j}(x)u_{j},
	\end{equation}
	where the weight function $\alpha(x)$=1, $\phi_{j}(x)$ is the Chebyshev polynomial of degree $\leq N$ and $u_{j}=u(x_{j})$.\\
	 The Chebyshev collocation points are defined by: 
	\begin{equation} \label{eq3} x_{j}=cos(\frac{j\pi}{N}),\quad j=0,....,N. \end{equation}
	A spectral differentiation matrix for the Chebyshev collocation points is created by interpolating a polynomial through the collocation points, i.e., the polynomial
	$$P_{N}(x_{k})= \sum_{j=0}^{N}\phi_{j}(x_{k})y_{j}.$$
	The derivatives values of the interpolating polynomial  \ref{eq2} at the Chebyshev collocation points  \ref{eq3} are:
	$$ P_{N}^{(l)}(x)=\sum_{j=0}^{N}\phi_{j}^{(l)}(x_{k})y_{j}.$$The differentiation matrix $D^{(l)}$ with entries 
	$$D^{(l)}_{k,j}=\phi_{j}^{(l)}(x_{k})$$ is explicitly determined in ~\cite {TrefethenF} and  ~\cite {Weiden}.\\
	For further information on convergence rates, collocation differentiation matrices, and the efficiency of the Chebyshev collocation method, refer to ~\cite{CANUTO,TrefethenF, Weiden,Gottlieb}
	
	\section{Numerical computations}
	
	\subsection {The Orr-Sommerfeld stability equation for plane Poiseuille flow}
	
	The stability of plane Poiseuille flow to disturbances of finite amplitude is affected by the characteristics of the higher-order modes of the Orr-Sommerfeld equation:
	 \begin{equation}\label{eq4}
		\frac{d^{4}u}{dx^{4}}-2\alpha^{2}\frac{d^{2}u}{dx^{2}}+\alpha^{4}u-i\alpha R \left[(U-\lambda)(\frac{d^{2}u}{dx^{2}}-\alpha^{2}u)-\frac{d^{2}U}{dx^{2}} u\right] =0, \quad -1<x<1
	\end{equation}
with boundary conditions $u(1)=u'(1)=0$, $u(-1)=u'(-1)=0$
where $U(x) = 1 - x^2$ is the base velocity, $\alpha$ is the wave number of the
disturbance, $R$ is the Reynolds number of the flow, and $\lambda$ is the temporal eigenvalue. \\
For a given Reynolds number,
the flow is stable if $Im(\lambda) < 0$ for all wave numbers $\alpha$ .

\begin{figure}[h]
	\caption{Orr-Sommerfeld eigenvalues for plane Poiseuille flow computed by Chebfun when $R = 10000$ and $ \alpha =1$}	
	\includegraphics[width=10cm]{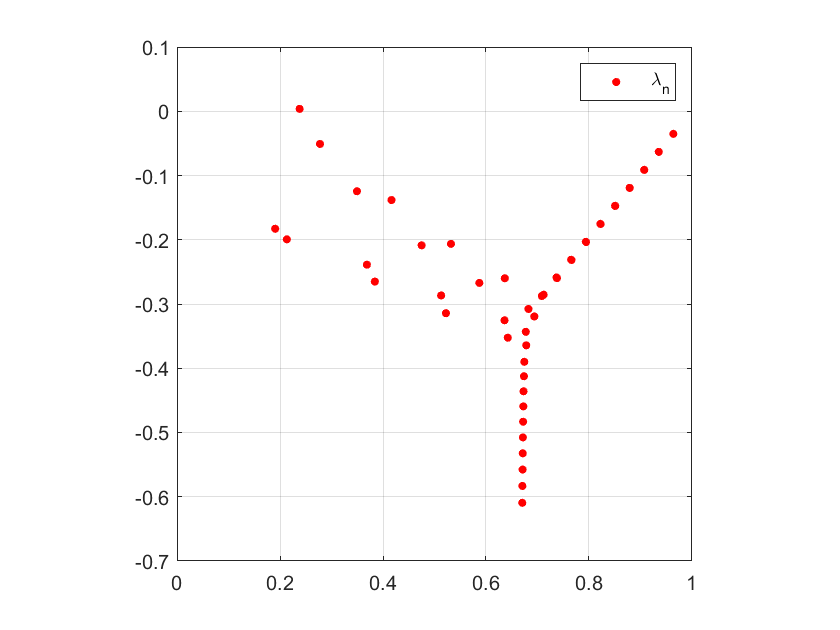}
	\label{fig:Fig1}
\end{figure}

\begin{figure}[h]
	\begin{center}
	\caption{Orr-Sommerfeld eigenfunctions for plane Poiseuille flow computed by Chebfun for n = 10, 20, 40, 50, $R = 10000$ and $ \alpha =1$}	
	\includegraphics[width=15cm]{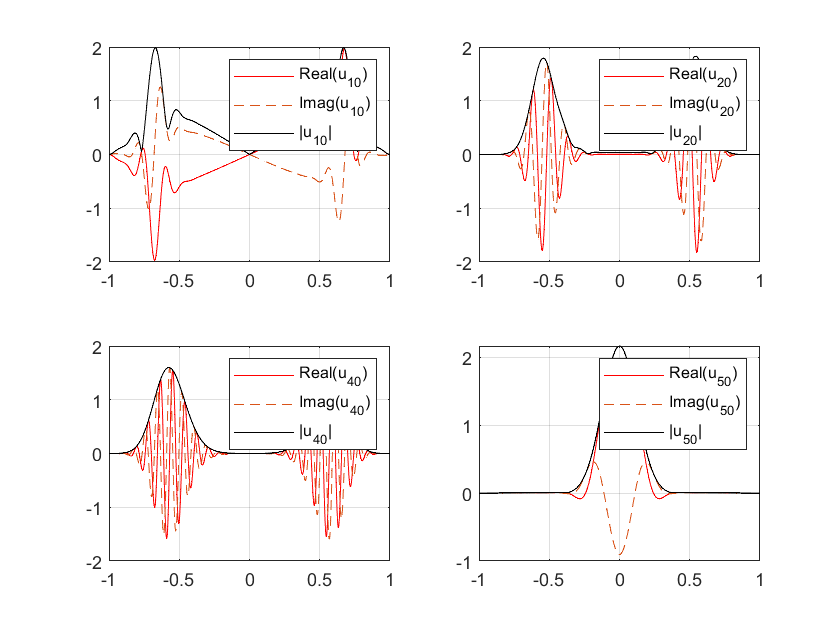}
	\label{fig:Fig2}
		\end{center}
\end{figure}

The Orr-Sommerfeld stability equation for plane Poiseuille flow has been solved by different methods, including the spectral Tau method ~\cite{Steven,Gardner,Fadden}. Despite its computational power, the Chebyshev tau method is prone to producing spurious eigenvalues, which can mislead by suggesting the presence of physical instabilities. We have solved  \ref{eq4} using the following Chebfun code:
\begin{lstlisting}[style=Matlab-editor]
      Re = 10000;                  
      alph = 1;                   
      A = chebop(-1,1);
      A.op = @(x,u) (diff(u,4)-2*alph^2*diff(u,2)+alph^4*u)/Re - ...
      2i*alph*u - 1i*alph*(1-x^2)*(diff(u,2)-alph^2*u);
      B = chebop(-1,1);
      B.op = @(x,u) (diff(u,2) - alph^2*u);
      A.lbc = [0; 0];
      A.rbc = [0; 0];
      [V,D] = eigs(A,B,50,'LR');
      e = diag(D)/(-i)
      [e,ii] = sort(e); V = V(:,ii)
\end{lstlisting} 
In Table \ref{tab:Table1}, the eigenvalues computed from Chebfun are compared with those obtained by Orzag in ~\cite{Steven}.\\
The numerical results in Table \ref{tab:Table1}, obtained using Chebfun algorithms, can significantly improve accuracy.\\
The first fifty eigenvalues are displayed in Figure \ref{fig:Fig1}. It is clear that Chebfun considerably improves the accuracy. These results demonstrate a significant improvement in convergence, with the obtained eigenvalues closely matching the exact values.\\
Figure \ref{fig:Fig2} shows the numerical computations of the Orr-Sommerfeld eigenfunctions for plane Poiseuille flow for n = 10, 20, 40, 50, $R = 10000$ and $ \alpha =1$.
\begin{table}[h]
	\caption{Computations of the first fifty eigenvalues $\lambda_{n}$ to problem \ref{eq4} and the relative error  with $R = 10000$ and $ \alpha =1$}
	\centering\begin{tabular}{c c c c  }
		\hline
		& & &  \\
		n\quad \quad  & $\lambda_{n}$ Current work & $\lambda_{n}$ ( ~\cite {Steven}) \quad \quad & Relative error $E_{n}$ \\
		\hline
		& & &  \\
1 & 0.237526453467190 + 0.003739621558176i   &	0.23752649 + 0.00373967i	& 0.002554066316713 e-07   \\
2 & 0.190059315832009 - 0.182821913447955i  & 0.1900592 - 0.1828219i	&  0.004421796894968 e-07 \\
3 & 0.277204423852887 - 0.050898788285925i  & 0.27720434 - 0.05089873i	&  0.003623362608429 e-07 \\
4 & 0.212725758914699 - 0.199360772538066i  &	& \\
5 & 0.349106782696747 - 0.124501988550457i  &	0.34910682 - 0.12450198i	&  0.001032547319289 e-07 \\
6 & 0.416351043011009 - 0.138226480657558i  &	0.41635102 - 0.13822652i	& 0.001038936208719  e-07\\
7 & 0.368498628885649 - 0.238824933625706i  &	0.36850 - 0.23882i	&  0.116610393163039 e-07\\
8 & 0.383987611550432 - 0.265106636966566i  &	0.38399 - 0.26511i	&   0.088399782260441 e-07\\
9 & 0.474900594285341 - 0.208731679514684i  &	0.474901 - 0.208731i	&  0.015256332571211 e-07\\
10 & 0.532045632398260 - 0.206464958879861i  &	0.532045 - 0.206466i	&   0.021344566260602 e-07\\
11 & 0.512917838578008 - 0.286625804918651i  &	0.51292 - 0.28663i	&   0.080315922225781 e-07\\
12 & 0.522227242060733 - 0.314297221673093i  &		&\\
13 & 0.587212152858069 - 0.267155024543561i &	0.58721 - 0.26716i	&  0.084033738591391 e-07\\
14 & 0.636726764020225 - 0.259891576417492i &	0.63672 - 0.25988i	&  0.194959976700758 e-07\\
15 & 0.635972403349023 - 0.325240246145666i &		&\\
16b& 0.642467128142601 - 0.352611661544509i &		&\\
17 & 0.682815324308629 - 0.307396236337561i &		&\\
18 & 0.677829693764668 - 0.343862562139100i &		&\\
19 & 0.694012915495241 - 0.319512596524431i &		&\\
20 & 0.708900943566356 - 0.287700499188017i &	0.70887 - 0.28765i	&  0.774181523230390 e-07\\
21 & 0.712380015769697 - 0.285532663566736i &   0.71232 - 0.28551i	& 0.835961539597376 e-07\\
22 & 0.678460650209095 - 0.364105997346648i &		&\\
23 & 0.674481311788510 - 0.389816181831927i &		&\\
24 & 0.737396932060859 - 0.258717338950187i &	0.73741 - 0.25872i	&  0.170652356699889 e-07\\
25 & 0.738110703108953 - 0.259683790147475i &	0.73812 - 0.25969i	& 0.142881964647236 e-07\\
26 & 0.674074031819750 - 0.412476248961745i &		&\\
27 & 0.765880870661141 - 0.231187092351485i &	0.76588 - 0.23119i	&0.037939443058216 e-07\\
28 & 0.766493381894095 - 0.231589331211106i &	0.76649 - 0.23159i	&0.043053997640574 e-07\\
29 & 0.673214263544665 - 0.435825255094054i &		&\\
30 & 0.672842907679867 - 0.459368746508057i &		&\\
31 & 0.794388205614881 - 0.203219928274886i &	0.794388 - 0.203221i	&  0.013308666078341 e-07\\
32 & 0.794819084348003 - 0.203528365908084i &	0.794818 - 0.203529i	& 0.015310111824788 e-07\\
33 & 0.672314283393095 - 0.483242312506917i &		&\\
34 & 0.822834837645231 - 0.175228699833616i &	0.8228350 - 0.1752287i	&0.001929840841948 e-07\\
35 & 0.823136863367242 - 0.175478000745007i   & 0.8231370 - 0.1754781i  &0.002006560203069 e-07\\
36 & 0.671999256182045 - 0.507679783156294i & &\\
37 & 0.671595417579113 - 0.532415153203374i& &\\
38 & 0.851245838857825 - 0.147233965973347i   & 0.8512458 - 0.1472339i& 0.000886303149438 e-07\\
39 & 0.851449336665746 - 0.147425594506901i    &0.8514494- 0.1474256i& 0.000735686572235 e-07\\
40 & 0.671338080555365 - 0.557660822424268i  & &  \\
41 & 0.879627300060967 - 0.119232862448449i  & 0.87962729 - 0.11923285i& 0.000180312665151 e-07\\
42 & 0.879755694612466 - 0.119370746639358i  & 0.87975570- 0.11937073i& 0.000196997956516 e-07\\
43 & 0.670961078894129 - 0.583266696359715i&&\\
44 & 0.670767188245712 - 0.609389827537272i&&\\
45 & 0.907983054287578 - 0.091222734983170i  &   0.90798305 - 0.09122274i& 0.000072317769954 e-07\\
46 & 0.908056340249451 - 0.091312862207931i   &  0.90805633 - 0.09131286i& 0.000700121296357  e-07\\
47 & 0.936316577176816 - 0.063201554173122i    & 0.93631654- 0.06320150i& 0.000114882280359 e-07\\
48 & 0.936351870150086 - 0.063251781457370i    & 0.93635178 - 0.06325157i& 0.002449395472099 e-07\\
49 & 0.964630915413379 - 0.035167277280544i   & 0.96463092 - 0.03516728i& 0.000055240575116 e-07\\
50 & 0.964642508235966 - 0.035186582968917i    &0.96464251 - 0.03518658i& 0.000035776469750 e-07\\
\hline
\end{tabular}
\label{tab:Table1}
\end{table}
	\subsection {Sixth-order eigenvalue problem}  
\begin{figure}[h]
	\caption{ Eigenvalues of a sixth-order eigenvalue problem \ref{eq5} computed by Chebyshev differentiation matrix  when the order of approximation has been N = 300}	
	\includegraphics[width=10cm]{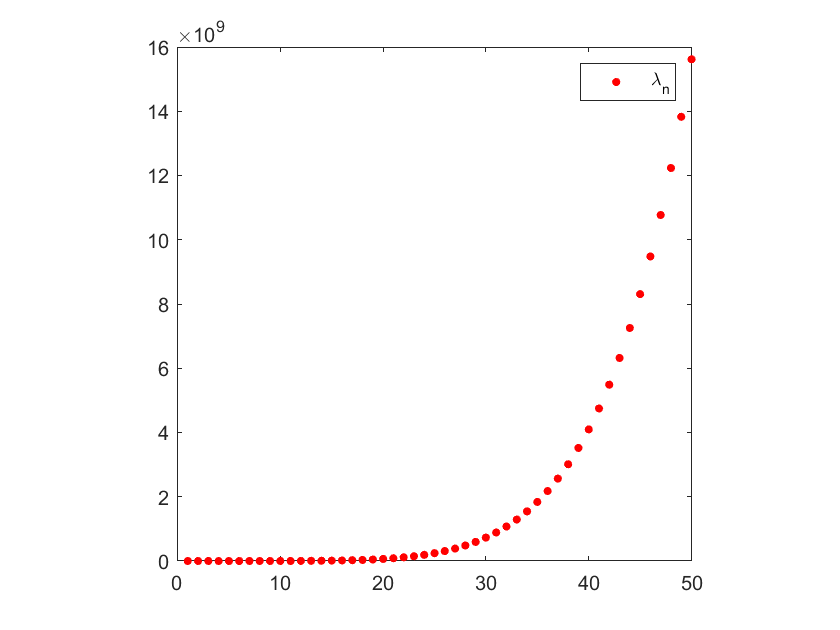}
	\label{fig:Fig3}
\end{figure}

\begin{figure}[h]
	\caption{Eigenfunctions of a sixth-order eigenvalue problem \ref{eq5} computed by Chebyshev differentiation matrix when the order of approximation has been N = 300}	
	\includegraphics[width=15cm]{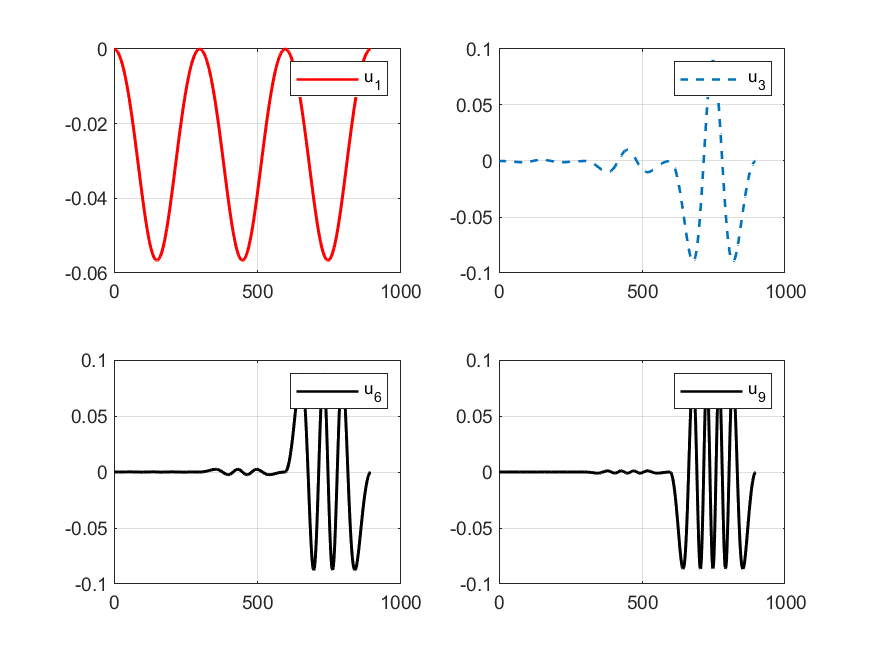}
	\label{fig:Fig4}
\end{figure}

	We consider a sixth-order eigenvalue problem of the form
	\begin{equation}\label{eq5}
		-u^{(vi)}=\lambda u, x\in[0,\pi]
	\end{equation}
with boundary conditions $u(0)=u''(0)=u^{(iv)}(0)=0$ and  $u(\pi)=u''(\pi)=u^{(iv)}(\pi)=0$.\\
This problem has exact eigenvalues $\lambda_{n}=(n+1)^{6}$, $n=0,1,2....$
We rewrite \ref{eq5} as a second
order differential system, namely
	\begin{equation*} 
	 \begin{array}{c c c  }
		u''(x)=& v(x)& x\in[0,\pi]\\
		v''(x) =& w(x)&\\
		w''(x)=& \lambda u(x)& \\
	\end{array} 
\end{equation*}
with boundary conditions $	u(0)=v(0)=w(0)=0$ and $ u(\pi)= v(\pi)=w(\pi)=0 $.
We construct the block matrices defined by
	\begin{equation*} \label{eq6}
	\left( \begin{array}{c c c}
		-\frac{4}{\pi^2}D_{N}^{(2)}&-I& O \\
		O & 	-\frac{4}{\pi^2}D_{N}^{(2)}&-I\\
		O& O& -\frac{4}{\pi^2}D_{N}^{(2)}
	\end{array} \right)
\end{equation*}
and 
	\begin{equation*} \label{eq7}
	B=\left( \begin{array}{c c c}
		O&O& O \\
		O & O&O	\\
		O& O&O 
	\end{array} \right)
\end{equation*}
The Chebyshev Collocation approach for solving \ref{eq5} involves the construction of 
	a $(N + 1) \times (N + 1)$  second order Chebyshev
	differentiation matrix $D_{N}^{(2)}$ associated with the nodes \ref{eq3},
	but shifted from the canonical Chebyshev interval $[-1,1]$ to $[0,\pi]$. \\
	The incorporation of the boundary conditions
 requires that the first and last rows of the matrix $D_{N}^{(2)}$ are removed,
	as well as its first and last columns (see ~\cite{FORNBERG}).\\
	The factor $\frac{4}{\pi^2}$ comes from the shift of interval $[-1,1]$ to  $[0,\pi]$.
	 $I$ and $O$ are, respectively, the identity and zero matrices of the same
	dimension as $D_{N}^{(2)}$.
	The following  MATLAB code has been used to solve \ref{eq5}:
	\begin{lstlisting}[style=Matlab-editor]
	N=300; 
	[x,D]=chebdif(N,2);
        D2=D(2:N-1,2:N-1,2); 
	I=eye(size(D2)); O=zeros(size(D2));
	k = 4/(pi^2); 
	A=[-k*D2 -I O; O -k*D2 -I; O O -k*D2];
	B=[O O O; O O O; I O O];  
	k = 10 ; 
	E=eigs(@(x)(A\(B*x)),k))
	\end{lstlisting}
The numerical results presented in Table \ref{tab:Table2} demonstrate the efficiency of the current technique. Figure \ref{fig:Fig3} displays the numerical computations of the eigenvalues. These results highlight the accuracy of the algorithms.\\
Note that this problem, as reported by Greenberg and Marletta ~\cite {Greenberg}, exhibits stiffness in at least part of the range. The approximated eigenvalues are compared with those obtained by ~\cite {Greenberg} in Table \ref{tab:Table2}. It is evident from the numerical results in Table \ref{tab:Table2} and Figure \ref{fig:Fig3} that there is a close agreement with the exact eigenvalues. Figure \ref{fig:Fig4} shows the numerical computations of the eigenfunctions associated to $\lambda_{1} $, $\lambda_{3} $, $\lambda_{6}$ and  $\lambda_{9} $.
\begin{table}[h]
		\caption{Computations of the first ten eigenvalues $\lambda_{n}$ to problem \ref{eq5}}
		\centering\begin{tabular}{c c c c }
			\hline
			& & &   \\
			n\quad \quad  & $\lambda_{n}$ Current work & Analytical eigenvalues $\lambda_{n}$ \quad \quad& $\lambda_{n}$( ~\cite {Greenberg})\quad \quad\\
			\hline
			& & &   \\
			0	& 0.999999999999436	&	1		&	0.9999998 \\
			1	&		63.9999999999705		&	64	& 63.999989	\\
			2	&	728.999999994709	&	729		&\\
			3		& 4095.99999999408		&	4096		&4095.9993	\\
			4	&		15625.0000025236	&	15625 		&	\\
			5		& 	46656.0000010131		& 	46656		&	46655.992	 	\\
			6		& 	117648.99948254			&  117649 		&		\\
			7	&		262143.999834618	& 262143.96		&	262143.96	\\
			8		& 	531441.006957261	& 531441	&			\\
			9		& 	999999.997801568		& 1000000	&	999999.83		\\
			& & & \\ 
			\hline
	\end{tabular}
	\label{tab:Table2}
\end{table}
	\subsection{Benard Stability Eigenvalue problem}  
	We consider the eigenvalue operator that arises in the analysis of Benard problem ~\cite {Greenberg} 
	\begin{equation}
		\label{eq8}
		(\frac{d^{2}}{dx^{2}}-a^{2})(\frac{d^{2}}{dx^{2}}-a^{2}-\nu)(\frac{d^{2}}{dx^{2}}-a^{2}-
		\dfrac{\nu}{P})u-a^{2}Ru= \lambda u , x\in[0,1]
	\end{equation}           
with boundary conditions $u(0)=u''(0)=u^{(iv)}(0)=0$ and  $u(1)=u''(1)=u^{(iv)}(1)=0$.\\
where $a$,$ P$, $\nu $ and $R$ are constants, and $\nu$ is regarded as the eigenparameter. \\
The eigenvalues  $\lambda_{n}$ , $n=0,1,2,...$ are functions of $\nu $, and the $\nu$-eigenvalues are given by
$$\nu^{\ast}_{n}=-\dfrac{(1+P)(n^{2}\pi^{2}+a^{2})}{2}\pm \sqrt{\dfrac{(1-P)^{2}(n^{2}\pi^{2}+a^{2})^{2}}{4}+\frac{a^{2}RP}{n^{2}\pi^{2}+a^{2}}},\quad n=0,1,2,...  $$
We solved the eigenproblem \ref{eq8} with $a = P = R = 1$ for several different values of $\nu $ using the Chebyshev collocation method. The eigenvalues are reported in Table \ref{tab:Table3}.
Note that $\nu^{\ast}_{1}=-(1+\pi^{2})$, $\nu^{\ast}_{2}=-(1+4\pi^{2})$, $\nu^{\ast}_{3}=-(1+\pi^{2})-(1+\pi^{2})^{-\frac{1}{2}}$, and  $\nu^{\ast}_{4}=-(1+4\pi^{2})-(1+4\pi^{2})^{-\frac{1}{2}}$.
	\begin{table}[h]
	\caption{Computations of the first two eigenvalues of Bénard problem \ref{eq8} for various $\nu^{\pm}_{n}$}
	\centering\begin{tabular}{c c c c c c c }
		\hline
		& & &  & &  & \\
		$\nu^{\pm}_{n}$\quad \quad  & $\lambda_{0}$  & $\lambda_{1}$ \quad \quad & $\lambda_{0}$ ~\cite {Greenberg} & $\lambda_{1}$~\cite {Greenberg} & $\lambda_{0}$ ~\cite {lesnic} & $\lambda_{1}$~\cite {lesnic} \\
		\hline
		& & & & &  & \\
		$-(1+\pi^{2})-(1+\pi^{2})^{-\frac{1}{2}}$&	-1.999999999379	& -34764.361138564&	$-1.10^{-7}$&34762.361&$-$ &$-$\\
		$-(1+\pi^{2})$ &	-1.00000000111486&	-35487.6927732063	&-1.000005&106725.18&$-1$& $ 4.81423\times 10^{10} $ \\
		$-(1+4\pi^{2})-(1+4\pi^{2})^{-\frac{1}{2}}$&	 -2.00000004492031&	-9631.62337424756&	$-3.10^{-5}$& 9629.62&$-$  &$-$\\		
		$-(1+4\pi^{2})$	&  -1.00000004976755 & -9530.18456147846&-1.0001 & 9528.17& $-1$ &$1,01881 \times10^{11}$\\
		\hline
	\end{tabular}
	\label{tab:Table3}
\end{table}

\begin{figure}[h]
	\caption{Relative error with respect to $N$ of the first thirty eigenvalues of the problem \ref{eq8} computed by Chebychev collocation method for $\nu^{\ast}_{4}=-(1+4\pi^{2})-(1+4\pi^{2})^{-\frac{1}{2}}$, $N_{1}=250$ and $N_{2}=500$. }	
	\includegraphics[width=15cm]{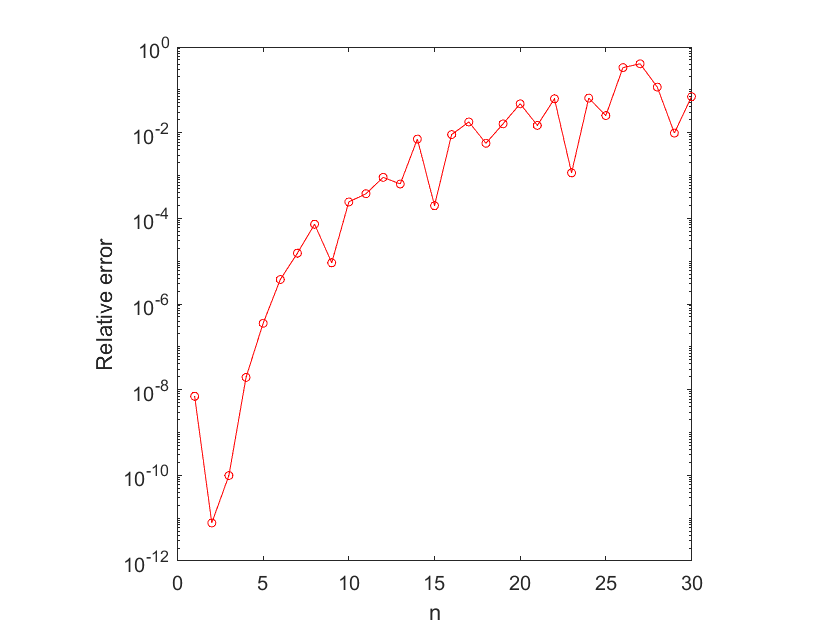}
	\label{fig:Fig9}
\end{figure}

Considering the approximations obtained for $\lambda_{0}$ which is approximately
equal to $-1$ for $\nu^{\ast}_{1}$ and $\nu^{\ast}_{2}$, it is clear that our results are very close to the findings reported in ~\cite {Greenberg}. It no longer happens for $\lambda_{1}$. Greenberg and Marletta ~\cite {Greenberg} mentioned that their code is not very accurate and Lesnic and Attili ~\cite {lesnic} did not approximate $\lambda_{0}$ for $\nu^{\ast}_{3}$ and $\nu^{\ast}_{4}$. However, based on the relative error concerning the order of approximation $N$ of the first thirty eigenvalues of the eigenproblem \ref{eq8} computed by Chebychev collocation method for $\nu^{\ast}_{4}=-(1+4\pi^{2})-(1+4\pi^{2})^{-\frac{1}{2}}$ and displayed in Figure \ref{fig:Fig9}, it suggests that $\lambda_{0}$ and $\lambda_{1}$ are computed with an accuracy of at least $10^{-11}$ and the first thirty with an accuracy of at least $10^{-1}$. We can assert that our eigenvalue approximations are more accurate than the results reported in ~\cite {Greenberg}  and ~\cite {lesnic}. 
	\subsection{The vibrations of a hinged wedge-shaped beam}
\begin{figure}[h]
	\caption{Eigenfunctions of the problem \ref{eq9} computed by Chebfun}	
	\includegraphics[width=15cm]{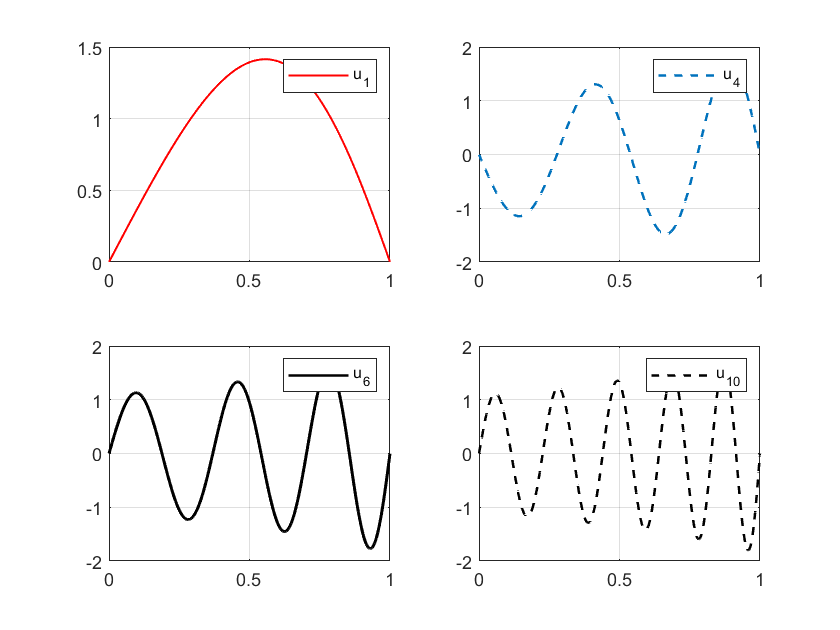}
	\label{fig:Fig5}
\end{figure}

The fourth order eigenproblem \begin{equation}
	\label{eq9}
	((1-Ax)^{3}u^{(2)})^{(2)}-\lambda(1-Ax)u=0, \quad\quad u(0)=u''(0)= u(1)=u''(1)=0. \end{equation}
The taper of the beam is described by the parameter $A$. This problem was considered in ~\cite {Banks} where the authors used the method of Prüfer Transformation to compute the first three eigenvalues. Tricomi has also analyzed this problem in ~\cite {Tricomi}.\\
 We apply Chebfun algorithms to compute the eigenvalues for $A = 0.5$ and compare our results with the values in ~\cite {Banks}. We report in Table \ref{tab:Table4} the latter eigenvalues compared with the ones provided
by Chebfun.
\begin{table}[h]
	\caption{The first ten eigenvalues to problem \ref{eq9} computed by Chebfun and ~\cite {Banks}.}
	\centering\begin{tabular}{c c c  }
		\hline
		& &   \\
		n\quad \quad  & $\lambda_{n}$ Current work & $\lambda_{n}$( ~\cite {Banks})\quad \quad\\
		\hline
		& &   \\
		1	& 50.7154258546	&	50.716228		 \\
		2	&		838.2127293477		&838.20834	\\
		3	&	4222.2443499295	& 4222.2394\\
		4		& 13305.8158141082		&	$-$	\\
		5	&		32431.9786551817	&$-$ 			\\
		6		& 67184.7630480295		& $-$	 	\\
		7		& 	124389.0502298803		& $-$ 			\\
		8	&		212110.6501871448	&$-$		\\
		9		& 339656.2267803520	& 		$-$	\\
		10		& 	517573.3580451152		& $-$		\\
		& &  \\ 
		\hline
	\end{tabular}
	\label{tab:Table4}
\end{table}
Some eigenfunctions to problem \ref{eq9} computed by Chebfun are displayed in
Figure \ref{fig:Fig5}. They satisfy the boundary conditions assumed in this problem.
\subsection{The longitudinal vibrations of a cantilevered beam}
\begin{figure}[h]
		\caption{Eigenfunctions of the problem \ref{eq10} computed by Chebfun}	
		\includegraphics[width=15cm]{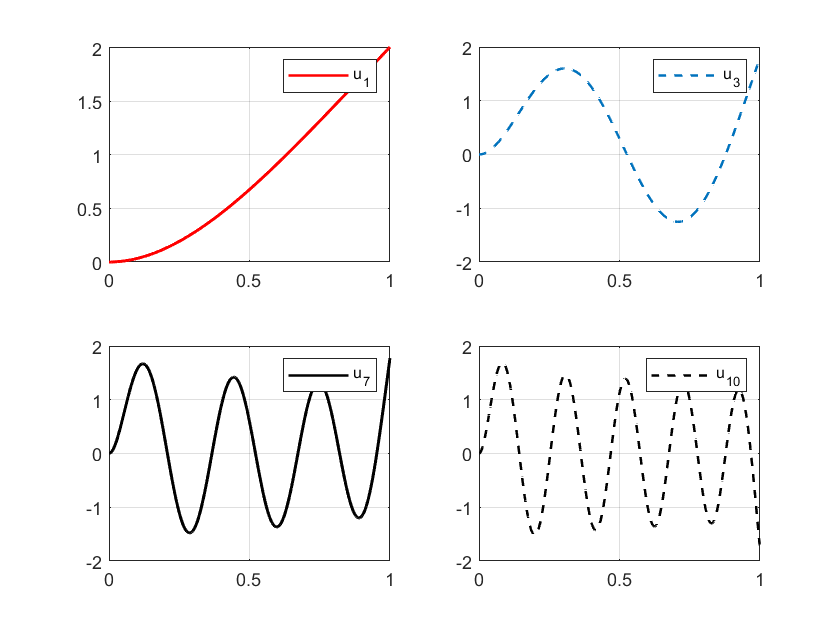}
		\label{fig:Fig6}
	\end{figure}

We consider the longitudinal vibrations of a cantilevered beam defined by the eigenvalue problem 
\begin{equation}	\label{eq10}
	-((1+x)u')'=\lambda(1+x)u,  \quad\quad u'(0)=0, u(1)=0. \end{equation}
This problem is considered in ~\cite{Collatz}, where the solution was obtained using finite differences.
To demonstrate the applicability of the Chebfun computational procedure for eigenvalues, as outlined in the previous problems, we compare our results with those in ~\cite {Collatz}. In Table \ref{tab:Table5}, we present the eigenvalues for problem \ref{eq10} and compare them with those provided by Chebfun.
	\begin{table}[h]
		\caption{The first ten eigenvalues to problem \ref{eq10} computed by Chebfun and ~\cite {Collatz}.}
		\centering\begin{tabular}{c c c  }
			\hline
			& &   \\
			n\quad \quad  & $\lambda_{n}$ Chebfun & $\lambda_{n}$(~\cite {Collatz})\quad \quad\\
			\hline
			& &   \\
			1	& 3.2184751303858	&	3·3731		 \\
			2	&		23.0597875591202		&23·19	\\
			3	&	62.5516753666466	& 62·68\\
			4		& 121.7733123600257		&	121·90	\\
			5	&		200.7318417346520	&200·86 			\\
			6		& 299.4287536221811	& 299·56	 	\\
			7		& 	417.8645088149235		& $-$ 			\\
			8	&	556.0392868523723	&$-$		\\
			9		& 713.9531692682501	& 		$-$	\\
			10		& 	891.6061973898490		& $-$		\\
			& &  \\ 
			\hline
		\end{tabular}
		\label{tab:Table5}
	\end{table}
	Some eigenfunctions for $n=1,3,7$ and $10 $ of the problem \ref{eq10} are displayed in
	Figure \ref{fig:Fig6}. 
	\subsection{The vibrations of a rotated beam}
\begin{figure}[h]
		\caption{Eigenvalues to problem \ref{eq11} computed by Chebfun}	
		\includegraphics[width=10cm]{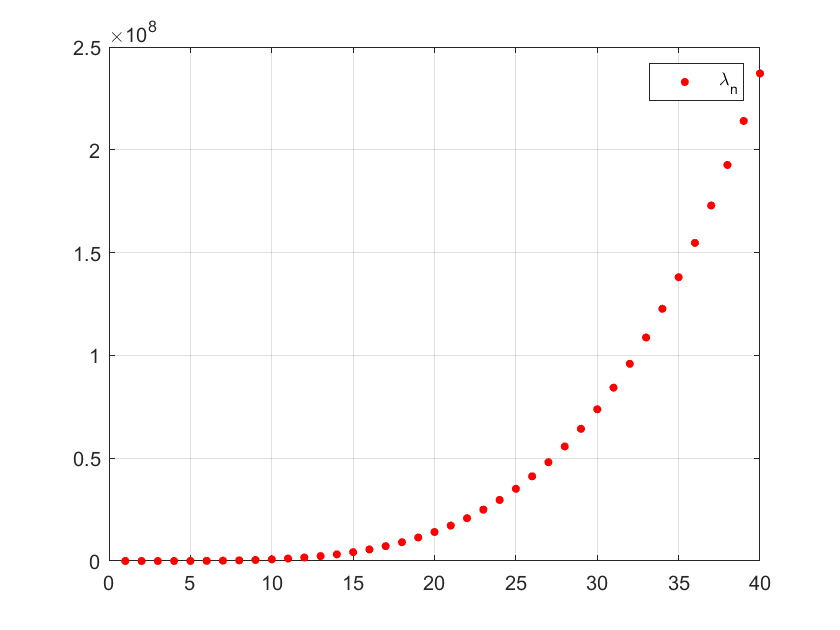}
		\label{fig:Fig7}
	\end{figure}

\begin{figure}[h]
	\caption{Eigenfunctions of the problem \ref{eq11} computed by Chebfun}	
	\includegraphics[width=15cm]{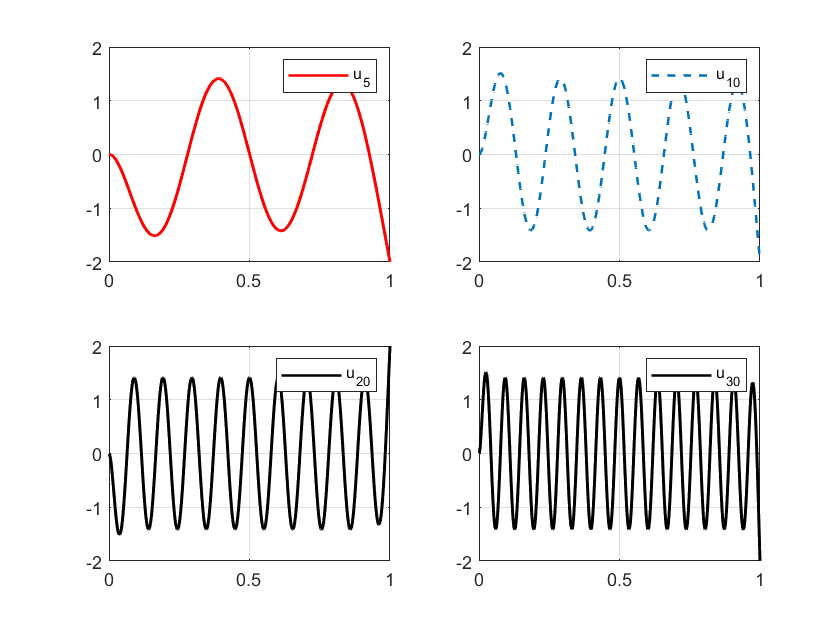}
	\label{fig:Fig8}
\end{figure}

We approximate the eigenvalues of a beam rotating about a fixed end in a plane that
contains the axis of the beam. The corresponding eigenvalue problem (see ~\cite{Boyce}) is \begin{equation}
			\label{eq11}
			u^{(iv)}-\frac{1}{2}\alpha^{2}((1-x)^{2} u')'=  \lambda u , x\in[0,1] \end{equation}
with boundary conditions $ u(0)=u'(0)=0, u''(1)=u'''(1)=0.$ 	Here $\alpha$ is proportional to the angular velocity of rotation. As an application of the decomposition of operators into sums (see ~\cite{Gould}), we consider the two separately resolvable eigenvalue problems 
		 \begin{equation}
		 	\label{eq12}
		 	u^{(iv)}=  \lambda u  , x\in[0,1]
		 \end{equation}
		and 
			\begin{equation}
			\label{eq13}
			-\frac{1}{2}\alpha^{2}((1-x)^{2} u')'=  \lambda u  , x\in[0,1]
		\end{equation}
		with boundary conditions $ u(0)=u'(0)=0, u''(1)=u'''(1)=0.$\\
		The eigenvalues $\lambda_{n}^{(1)}$ and eigenfunctions $u_{n}^{(1)}$ of the eigenvalue problem  $\ref{eq12}$  are given by
		$$\lambda_{n}^{(1)}=(\beta_{n})^{4},$$
		$$u_{n}^{(1)}(x)=cosh (\beta_{n} x) -cos (\beta_{n} x)-\frac{cosh (\beta_{n})  +cos (\beta_{n}) }{sinh (\beta_{n})  +sin(\beta_{n}) }(sinh (\beta_{n}x)  -sin(\beta_{n}x) ) , n=1,2,3,....$$
		where $ \beta_{n}$ are the positive roots, in order of magnitude, of the equation
		$$cosh(\beta) cos (\beta)+1=0$$
		The eigenvalues $\lambda_{n}^{(2)}$ and eigenfunctions $u_{n}^{(2)}$ of the eigenvalue problem  $\ref{eq13}$   are given by
		$$\lambda_{n}^{(2)}=\alpha^{2}n(2n-1),$$
		$$u_{n}^{(2)}(x)=(-1)^{n}(4n-1)^{\frac{1}{2}}p_{2n-1}(x), n=1,2,3,....$$
		where $p_{2n-1}$ is the Legendre polynomial of degree $2n-1$.
		Since Chebfun can cope with various boundary conditions, we use the following code to solve the eigenvalue problem $\ref{eq11}$ for $\alpha^2=5$.
		\begin{lstlisting}[style=Matlab-editor]
	dom=[0,1];
	L = chebop(dom);
	x = chebfun('x',dom);
	L.op = @(x,u) diff(u,4)+5*x*diff(u)-(5/2)*(1-((x)^2))*diff(u,2);
	L.lbc = @(u)[u; diff(u,1)]; % fixed b. c.
	L.rbc = @(u)[diff(u,2); diff(u,3)];% free b. c.
	[U,D]=eigs(L,10)
	\end{lstlisting}
		In Table \ref{tab:Table6} and Figure \ref{fig:Fig7}, the first ten eigenvalues computed by chebfun and Bazley ~\cite {Bazley} are reported. The eigenfunctions are displayed in Figure  \ref{fig:Fig8}.
	The results indicate that the current algorithm performs slightly better than the variational method in ~\cite {Bazley}. 
			\begin{table}[h]
			\caption{The first ten eigenvalues to problem \ref{eq11} computed by Chebfun and ~\cite {Bazley}.}
			\centering\begin{tabular}{c c c  }
				\hline
				& &   \\
				n\quad \quad  & $\lambda_{n}$ Chebfun & Lower bounds $\lambda_{n}$(~\cite {Bazley})\quad \quad\\
				\hline
				& &   \\
				1	& 18.2525379005872	&	18.287		 \\
				2	&		518.105590004563		&517.41	\\
				3	&	3893.03175579422	& 3891.4\\
				4		& 14797.3127489787		&	14784	\\
				5	&		40247.3614834009	&$-$ 			\\
				6	&		89596.0904823395	&$-$			\\
				7		& 	174530.478417044& $-$ 	\\
				8		& 	309081.072633057		& $-$ 			\\
				9	&	509607.922852811	&$-$		\\
				10		& 794818.779148833	& 		$-$	\\
				& &  \\ 
				\hline
			\end{tabular}
			\label{tab:Table6}
		\end{table}
	\section{Conclusion}
	In this paper, we explore the application of Chebfun and Chebyshev collocation methods to solve high-order eigenvalue problems with a particular focus on their performance under various boundary conditions.
	One of the most compelling findings of our investigation is Chebfun's versatility in handling high-order boundary conditions for fourth-order eigenvalue problems. This capability is a notable advantage, as it allows researchers to apply Chebfun directly to various problems without requiring extensive modifications or adaptations. The reliability of the outcomes produced by Chebfun under these conditions is a testament to its robustness as a computational tool.
	
	However, the situation becomes less favorable when addressing sixth-order eigenvalue problems. The direct application of Chebfun in these cases poses challenges, particularly due to the emergence of poorly conditioned matrices. For sixth-order problems, especially those with hinged boundary conditions, reducing the problem to a second-order system is more effective before applying the Chebyshev collocation methods.

	\textbf{Conflicts of Interest:}\\ 
	The author declares that there are no conflicts of interest regarding the publication of this paper.
	
	\textbf{Data availability:}\\  The author confirms that the data presented in this study are available in supplementary materials.
	
	\textbf{Funding:}  This research received no external funding.
	\bibliographystyle{unsrtnat}
	\bibliography{references}
\end{document}